\documentclass[12pt,twoside]{amsart}
\usepackage{amsthm,amsfonts,amssymb,amscd,euscript}

\usepackage[matrix,arrow]{xy}

\author{Florin Ambro} 
\address{RIMS, Kyoto University\\
Kyoto 606-8502, Japan.}
\email{ambro@kurims.kyoto-u.ac.jp}

\newcommand{\C}{{\mathbb C}}

\newcommand{\Q}{{\mathbb Q}}
\newcommand{\Z}{{\mathbb Z}}

\newcommand{\R}{{\mathbb R}}

\newcommand{\emb}{\operatorname{emb}}

\newcommand{\ind}{\operatorname{index}}

\newcommand{\mld}{\operatorname{mld}}
\newcommand{\Mld}{\operatorname{Mld}}
\newcommand{\mult}{\operatorname{mult}}
\newcommand{\orb}{\operatorname{orb}}
\newcommand{\relint}{\operatorname{relint}}

\theoremstyle{plain}
\newtheorem{thm}{Theorem}[section]
\newtheorem{theom}{Theorem}
\newtheorem{lem}[thm]{Lemma}
\newtheorem{cor}[thm]{Corollary}

\newtheorem{prop}[thm]{Proposition}

\theoremstyle{definition}
\newtheorem{defn}[thm]{Definition}

\newtheorem{rem}[thm]{Remark}

\theoremstyle{remark}

\begin{document}

\bibliographystyle{amsalpha+}
\title[The set of toric minimal log discrepancies]
{The set of toric minimal log discrepancies}
%\date{June 23, 2004}

\begin{abstract} We describe the set of minimal log 
discrepancies of toric log varieties, and study its 
accumulation points.
\end{abstract}

\maketitle

%%%%%%%%%%%%%%%%%%%%%%%%%%%%%%%%%%%%%%%%%%%%%%%%%%%%%%%
%\setcounter{section}{-1}
%%%%%%%%%%%%%%%%%%%%%%%%%%%%%%%%%%%%%%%%%%%%%%%%%%%%%%%
%%% Document name: toricmld.tex
%%% Last modified: Wed Jan 11 18:24:11 JST 2006
%%%%%%%%%%%%%%%%%%%%%%%%%%%%%%%%%%%%%%%%%%%%%%%%%%%%%%%

%%%%%%%%%%%%%%%%%%%%%%%%%%%%%%%%%%
%%%%%%%%%%%%%%%%%%%%%%%%%%%%%%%%%%

\section*{Introduction}

%%%%%%%%%%%%%%%%%%%%%%%%%%%%%%%%%%%
%%%%%%%%%%%%%%%%%%%%%%%%%%%%%%%%%%%

\footnotetext[1]{This work is supported by a 
21st Century COE Kyoto Mathematics Fellowship,
and a JSPS Grant-in-Aid No 17740011.
}
\footnotetext[2]{1991 Mathematics Subject Classification. 
Primary: 14B05. Secondary: 14M25.}

Minimal log discrepancies are invariants of 
singularities of log varieties. A log variety 
$(X,B)$ is a normal variety $X$ endowed 
with an effective Weil $\R$-divisor $B$, having
at most log canonical singularities. For any 
Grothendieck point $\eta\in X$, the minimal 
log discrepancy of $(X,B)$ at $\eta$, denoted by
$a(\eta;X,B)$, is a non-negative real number.
For example, $a(\eta;X,B)=1-\mult_\eta(B)$ for 
every codimension one point $\eta\in X$. 
For higher codimensional points, minimal log 
discrepancies can be computed on a suitable 
resolution of $X$.

Let $A\subset [0,1]$ be a set containing $1$ and
let $d$ be a positive integer.
Denote by $\Mld_d(A)$ the set of minimal log 
discrepancies $a(\eta;X,B)$, where $\eta\in X$ is 
a Grothendieck point of codimension $d$, and $(X,B)$ 
is a log variety whose minimal log discrepancies 
in codimension one belong to $A$.
For example, $\Mld_1(A)=A$.
In connection to the termination of a sequence 
of log flips (see~\cite{Shokurov88, Termination}), 
Shokurov conjectured that if $A$ satisfies the
ascending chain condition, so does $\Mld_d(A)$. 
Furthermore, under certain
assumptions, the accumulation points of $\Mld_d(A)$ 
should correspond to minimal log discrepancies of smaller
codimensional points. This is known to hold
for $d=2$ (Shokurov~\cite{Shokurov93}, 
Alexeev~\cite{Alexeev93}) and for any $d$ in the case
of toric varieties without boundary 
(Borisov~\cite{Borisov97}). 
The purpose of this note is to extend Borisov's 
result to the case of toric log varieties. Given
the explicit nature of the toric case, we hope 
this will provide the reader with interesting examples.

In order to state the main result, define $\Mld^{tor}_d(A)
\subset\Mld_d(A)$ as above, except that we further require 
that $X$ is a toric variety and $B$ is torus invariant.
Note that $\Mld^{tor}_1(A)=A$.

\begin{theom}\label{main} The following properties hold
for $d\ge 2$:
\begin{itemize}
\item[(1)] We have
$$
\Mld^{tor}_d(A)
=\{\sum_{i=1}^s x_i a_i 
\left| \begin{array}{l}
2\le s\le d \\
(x_1,\ldots,x_s)\in \Q^s\cap (0,1]^s, (a_1,\ldots,a_s)\in A^s\\
\ind(x_i)\vert \ind(x_1,\ldots,\hat{x_i},\ldots,x_s),
\ \forall 1\le i\le s \\
\sum_{i=1}^s (1+(m-1)x_i-\lceil mx_i\rceil)a_i\ge 0\ \forall
m\in \Z
\end{array}\right
\},
$$ 
where for a rational point $x\in \Q^n$, we denote by
$\ind(x)$ the smallest positive integer $q$ 
such that $qx\in \Z^n$.
\item[(2)] If $A$ satisfies the ascending chain 
condition, then so does $\Mld^{tor}_d(A)$.
\item[(3)] Assume that $A$ has no nonzero accumulation
points. Then the set of accumulation points of 
$\Mld^{tor}_d(A)$ is included in
$$
\{0\}\cup \bigcup_{1\le d'\le d-1}\Mld^{tor}_{d'}
(\{\frac{1}{n};n\ge 1\}\cdot A).
$$ 
Equality holds if $d=2$, or if
$\{\frac{1}{n};n\ge 1\}\cdot A \subseteq A$.
\end{itemize}
\end{theom}

 We use the same methods as Borisov~\cite{Borisov97, 
Borisov99}. The explicit description in (1) is 
straightforward, whereas the accumulation behaviour 
in (2) and (3) relies on a result of 
Lawrence~\cite{Lawrence91} stating that the set of closed
subgroups of a real torus, which do not intersect a given 
open subset, has finitely many maximal elements with 
respect to inclusion.

 Finally, we should point out that $\Mld^{tor}_d(A)$ 
is strictly smaller than $\Mld_d(A)$ in general. For example,
even the set of accumulation points of $\Mld_2(A)$ (see
Shokurov~\cite{Shokurov93} for an explicit description) is
larger than $\{0\}\cup \{\frac{1}{n};n\ge 1\}\cdot A$, the
set of accumulation points of $\Mld^{tor}_2(A)$.

%%%%%%%%%%%%%%%%%%%%%%%%%%%%%%%%%%
%%%%%%%%%%%%%%%%%%%%%%%%%%%%%%%%%%

\section{Toric log varieties}

%%%%%%%%%%%%%%%%%%%%%%%%%%%%%%%%%%%
%%%%%%%%%%%%%%%%%%%%%%%%%%%%%%%%%%%

In this section we recall the definition of minimal log 
discrepancies and their explicit description in the toric 
case. The reader may consult~\cite{Ambro97} for more details.

A {\em log variety} $(X,B)$ consists of an algebraic variety $X$, 
defined over an algebraically closed field of characteristic 
zero, endowed with a finite combination $B=\sum_i b_i B_i$ 
of Weil prime divisors with real coefficients, such that $K_X+B$
is $\R$-Cartier. Here $K_X$ is the canonical divisor of $X$,
computed as the Weil divisor of zeros and poles $(\omega)_X$ of 
a top rational form $\omega\in \wedge^{\dim(X)}\Omega^1_{\C(X)/\C}$;
it is uniquely defined up to linear equivalence. The $\R$-Cartier
property of $K_X+B$ means that locally on $X$, there exists 
finitely many non-zero
rational functions $a_\alpha\in \C(X)^\times$ and $r_\alpha\in 
\R$ such that $K_X+B=\sum_\alpha r_\alpha(a_\alpha)$.

Let $\mu\colon X'\to X$ be a proper birational morphism from a
normal variety $X'$ and let $E\subset X'$ be a prime divisor.
Let $\omega$ be a top rational form on $X$, defining $K_X$,
and let $K_{X'}$ be the canonical divisor defined by 
$\mu^*\omega$. The real number
$$
a(E;X,B)=1+\mult_E(K_{X'}-\mu^*(K_X+B))
$$
is called the {\em log discrepancy of $(X,B)$ at $E$}. 
For a Grothendieck point $\eta\in X$, the {\em minimal
log discrepancy of $(X,B)$ at $\eta$} is defined as
$$
a(\eta;X,B)=\inf_{\mu(E)=\bar{\eta}}a(E;X,B),
$$
where the infimum is taken after all prime divisors $E$
on proper birational maps $\mu\colon X'\to X$. This infimum 
is either $-\infty$, or a non-negative real number. In the
latter case, $(X,B)$ is said to have {\em log canonical 
singularities at} $\eta$ and the invariant is computed as
follows. By Hironaka, there exists a proper birational 
morphism $\mu\colon X'\to X$ such that $X'$ is nonsingular,
$\mu^{-1}(\bar{\eta})$
is a divisor on $X'$, and there exists a simple normal
crossings divisor $\sum_i E_i$ on $X'$ which supports
both $\mu^{-1}(\bar{\eta})$ and $K_{X'}-\mu^*(K_X+B)$.
Then 
$$
a(\eta;X,B)=\min_{\mu(E_i)=\bar{\eta}}a(E_i;X,B).
$$

Next we specialize these notions to the toric case. 
We employ standard terminology on toric varieties, cf. 
Oda~\cite{Oda}.
A {\em toric log variety} is a log variety $(X,B)$ 
such that $X$ is a toric variety and $B$ is torus
invariant. Thus there exists a 
fan $\Delta$ in a lattice $N$ such that
$X=T_N\emb(\Delta)$ and $B=\sum_i b_i V(e_i)$, where
$\{e_i\}_i$ is the set of primitive lattice points 
on the one-dimensional cones of $\Delta$ and $V(e_i)\subset X$
is the torus invariant prime Weil divisor corresponding to $e_i$.
The canonical divisor is $K_X=\sum_i -V(e_i)$, and the $\R$-Cartier 
property of $K_X+B$ means that there exists a function
$\psi\colon \vert \Delta\vert\to \R$ such that $\psi(e_i)=1-b_i$
for every $i$, and $\psi\vert\sigma$ is linear for every cone 
$\sigma\in \Delta$. We may assume that $(X,B)$ has log
canonical singularities, which is equivalent to 
$\psi\ge 0$ or $b_i\in [0,1]$ for all $i$.

Let $e\in N^{prim}\cap \vert\Delta\vert$ be a non-zero primitive vector.
The barycentric subdivision with respect to $e$ defines a subdivision 
$\Delta_e\prec\Delta$ and the exceptional locus of the birational morphism 
$T_N\emb(\Delta_e)\to T_N\emb(\Delta)$ is a prime divisor denoted $E_e$. 
It is easy to see that 
$$
a(E_e;X,B)=\psi(e).
$$
Due to this property, $\psi$ is called the {\em log discrepancy 
function of} $(X,B)$.

Minimal log discrepancies of toric log varieties are computed
as follows. These are local invariants, so we 
only consider affine varieties. Thus $\Delta$ consists of the
faces of some strongly convex rational polyhedral cone
$\sigma\subset N_\R$ and we denote $X=T_N\emb(\sigma)$. Assume 
first that $0\in X$ is a torus invariant closed point
(it is unique since $X$ is affine). Using the existence of good
resolutions in the toric category, it is easy to see that
$$
a(0;X,B)=\min(\psi\vert_{N\cap \relint(\sigma)}).
$$
For the general case, let $\eta\in X$ be a Grothendieck point. 
There exists a unique face $\tau\prec \sigma$ such that 
$\eta\in \orb(\tau)$. Let $c$ and $d$ be the codimension 
of $\orb(\tau)$ and $\eta$ in $X$, respectively. The induced 
affine toric log variety 
$$
(X',B')=
(T_{N\cap(\tau-\tau)}\emb(\tau),\sum_{e\in\tau(1)}\mult_{V(e)}(B)V(e))
$$
has a unique torus invariant closed point $0'$, and we obtain
$$
a(\eta;X,B)=\mld(0';X',B')+d-c.
$$

%%%%%%%%%%%%%%%%%%%%%%%%%%%%%%%%%%
%%%%%%%%%%%%%%%%%%%%%%%%%%%%%%%%%%

\section{The set of toric minimal log discrepancies}

%%%%%%%%%%%%%%%%%%%%%%%%%%%%%%%%%%%
%%%%%%%%%%%%%%%%%%%%%%%%%%%%%%%%%%%

Let $A\subseteq [0,1]$ be a set containing $1$.

\begin{defn} For an integer $d\ge 1$, let $\Mld^{tor}_d(A)$ 
be the set of minimal log discrepancies $a(\eta;X,B)$, 
where $\eta\in X$ is a Grothendieck point of codimension 
$d$ and $(X,B)$ is a toric log variety whose minimal log 
discrepancies in codimension one belong to $A$. 
\end{defn}

It is easy to see that $\Mld^{tor}_1(A)=A$.

\begin{defn}
For an integer $d\ge 2$, define $V_d(A)$ to be the
set of pairs $(x,a)\in (0,1]^d\times A^d$ satisfying 
the following properties:
\begin{itemize}
\item[(i)] $x\in \Q^d$.
\item[(ii)] $\ind(x_i)\vert 
\ind(x_1,\ldots,\hat{x_i},\ldots,x_d)$ for $1\le i\le d$. 
\item[(iii)] $\sum_{i=1}^d(1+(m-1)x_i-\lceil mx_i\rceil)
a_i\ge 0$ for all $m\in \Z$.
\end{itemize}
For $x\in \Q^n$, $\ind(x)$ denotes the smallest positive 
integer $q$ such that $qx\in \Z^n$.
\end{defn}  

Note that property (ii) means that 
$(1,0,\ldots,0),\ldots,(0,\ldots,0,1)$ are primitive vectors 
in the lattice $\Z^d+\Z x$. Also, it is enough to verify 
property (iii) for the finitely many integers 
$1\le m\le \ind(x)-1$. For $(x,a)\in V_d(A)$ we denote 
$$
\langle x,a\rangle=\sum_{i=1}^d x_i a_i.
$$

\begin{prop}\label{gutaiteki} 
For $d\ge 2$, we have
$$
\Mld^{tor}_d(A)=\bigcup_{2\le s\le d}
\{\langle x,a\rangle; (x,a)\in V_s(A)\}.
$$
\end{prop}

\begin{proof} (1) We first show that the right hand side
is included in the left hand side. Fix $(x,a)\in V_s(A)$
for some $2\le s\le d$. 

If $s=d$, let $N=\Z^d+\Z x$ and let $\sigma$ be the standard 
positive cone in $\R^d$, spanned by standard basis 
$e_1,\ldots,e_d$ of $\Z^d$. Let $0\in T_N\emb(\sigma)$ be
the invariant closed point corresponding to $\sigma$. Then
the affine toric log variety
$$
(T_N\emb(\sigma),\sum_{i=1}^{d}(1-a_i)V(e_i))
$$
has minimal log discrepancy at $0$ equal to $\langle x,a\rangle$. 
Indeed, the log discrepancy function $\psi=\sum_{i=1}^d a_ie_i^*$ 
attains its minimum at $x$, and $\psi(x)=\langle x,a\rangle$. 
Therefore $\langle x,a\rangle\in \Mld^{tor}_d(A)$.

Assume now that $2\le s\le d-1$. Let $e_1,\ldots,e_d$ 
be the standard basis of $\Z^d$, let $e_{d+1}=(d-s)
e_1+e_2-\sum_{i=s+1}^d e_i$, let $v=\sum_{i=1}^s x_i 
e_i$ and let $N=\Z^d+\Z v$. Let $\sigma$ be
the cone in $\R^d$ generated by $e_1,\ldots,e_{d+1}$.
Set $a_i=a_1$ for $s+1\le i\le d$ and $a_{d+1}=a_2$. 
Then 
$$
0\in (T_N\emb(\sigma),\sum_{i=1}^{d+1}(1-a_i)V(e_i))
$$
is a $d$-dimensional germ of a toric log variety with 
minimal log discrepancy equal to $\langle x,a\rangle$.

Indeed, note first that the log variety is well defined
since $a_2=(d-s)a_1+a_2-\sum_{i=s+1}^d a_1$;
the log discrepancy function is $\psi=\sum_{i=1}^d a_i e_i^*$.
There exists $e=\sum_{i=1}^{d+1}y_i e_i\in N\cap 
\relint(\sigma)$ where the log discrepancy function $\psi$ 
attains its minimum. We may assume $y_i\in [0,1]$ for every $i$.
If $y_{d+1}\notin \Z$, then $y_{s+1}=\cdots=y_d=y_{d+1}$,
hence $e=\sum_{i=1}^s y_i e_i$. Therefore $\psi(e)\ge 
\psi(v)$. If $y_{d+1}\in \Z$, then $\sum_{i=1}^s y_i e_i
\in N\cap \relint(\sigma)$, hence  
$\psi(e)\ge \psi(\sum_{i=1}^s y_i e_i)\ge \psi(v)$.
We conclude that $\psi$ attains its minimum at $v$.
Therefore $\langle x,a\rangle=\psi(v)\in \Mld^{tor}_d(A)$.

(2) Let $(X,B)$ be a toric log variety with codimension one 
log discrepancies in $A$ and let $\eta\in X$ be a Grothendieck 
point of codimension $d$. We want to show that $a(\eta;X,B)$
belongs to the set on the right hand side.

There exists a unique cone $\sigma$ in the fan defining $X$ 
such that $\eta\in \orb(\sigma)$. Let $c$ be the codimension 
of $\orb(\sigma)$ in $X$. Then $a(\eta;X,B)$ coincides with 
the minimal log discrepancy of the toric log variety
$$
(T_{N\cap (\sigma-\sigma)}\emb(\sigma),\sum_{e\in
\sigma(1)}\mult_{V(e)}(B)V(e))\times {\mathbb C}^{d-c}.
$$
in the invariant closed point $0$. Therefore we may assume
that $X$ is affine and $\eta$ is a torus invariant closed
point $0$. 

We have $X=T_N\emb(\sigma)$, with $\dim N=d$, 
$B=\sum_{i\in I} (1-a_i)V(e_i)$ with $a_i\in A$ for every $i$.
The log discrepancy function $\psi\in \sigma^\vee$
of $(X,B)$ satisfies $\psi(e_i)=a_i$, and we have
$$
\mld(0;X,B)=\min(\psi\vert_{N\cap \relint(\sigma)}).
$$
There exists $e\in N\cap \relint(\sigma)$ such that
$\mld(0;X,B)=\psi(e)$. It is easy to see that there exists
a subset $\{1,\ldots, s\}\subseteq I$, with $2\le s\le d$, such 
that $e_1,\ldots,e_s$ are linearly independent and $e$ belongs 
to the relative interior of the cone spanned by 
$e_1,\ldots,e_s$. Let $e=\sum_{i=1}^s x_ie_i$, and denote 
$x=(x_1,\ldots,x_s)\in (0,1]^d, a=(a_1,\ldots,a_s)\in A^d$. 
It is clear that $\mld(0;X,B)=\langle x,a \rangle$, and we claim 
that $(x,a)\in V_s(A)$. 

Indeed, it is clear that $x\in \Q^s$. Since $e_i$ is a primitive 
lattice point of $N$, it is also a primitive in the sublattice 
$\sum_{i=1}^s \Z e_i+\Z e$, which is equivalent to 
$\ind(x_i)\vert \ind(x_1,\ldots,\hat{x_i},\ldots,x_s)$ for every 
$1\le i\le s$. Finally, let $m\in \Z$. We have
$\sum_{i=1}^s(1+mx_i-\lceil mx_i\rceil)e_i\in N\cap 
\relint(\sigma)$, hence
$\psi(\sum_{i=1}^s(1+mx_i-\lceil mx_i\rceil)e_i)\ge \psi(e)$.
Equivalently, $\sum_{i=1}^s(1+(m-1)x_i-\lceil mx_i\rceil)
a_i\ge 0$. Therefore $(x,a)\in V_s(A)$. 
\end{proof}

%%%%%%%%%%%%%%%%%%%%%%%%%%%%%%%%%%
%%%%%%%%%%%%%%%%%%%%%%%%%%%%%%%%%%

\section{The set $\tilde{V}_d(A)$}

%%%%%%%%%%%%%%%%%%%%%%%%%%%%%%%%%%%
%%%%%%%%%%%%%%%%%%%%%%%%%%%%%%%%%%%

By Proposition~\ref{gutaiteki}, the limiting behaviour 
of toric of minimal log discrepancies is controlled by 
the limiting behaviour of the sets $V_d(A)$. The rationality 
properties (i) and (ii) defining $V_d(A)$ do not behave 
well with respect to limits, and for this reason we enlarge
$V_d(A)$ to a new set $\tilde{V}_d(A)$, defined only by 
property (iii), which turns out to have good inductive 
properties and limiting behaviour.

\begin{defn}
Let $A\subseteq [0,1]$ be a subset containing $1$.
Define
$$
\tilde{V}_d(A)=\{(x,a)\in (0,1]^d\times A^d; 
\sum_{i=1}^d(1+(m-1)x_i-\lceil mx_i\rceil)a_i\ge 0,
\forall m\in \Z\}.
$$
\end{defn}

Equivalently, $\tilde{V}_d(A)$ is the set of pairs
$(x,a)\in (0,1]^d\times A^d$ such that the group
$\Z^d+\Z x$ does not intersect the set
$\{y\in (0,1]^d; \langle y-x,a\rangle<0\}$. 
As before, we denote $\langle x,a\rangle=\sum_{i=1}^d x_i a_i$.

\begin{lem}\label{onedim} 
The following equality holds
$$
\tilde{V}_1(A)=((0,1]\times \{0\})\cup(\{\frac{1}{n}; n\ge 1\}
\times A),
$$
where the first term is missing if $0\notin A$. In 
particular,
$$
\{\langle x,a\rangle; (x,a)\in \tilde{V}_1(A)\}=
\bigcup_{n=1}^\infty \frac{1}{n}\cdot A.
$$
\end{lem}

\begin{proof} Let $x\in (0,1]$ such that
$1+(m-1)x-\lceil mx\rceil\ge 0$ for every integer $m$.
Equivalently, we have
$$
\sup_{m\in \Z}(\lceil mx\rceil-mx) \le 1-x.
$$
Assume by contradiction that $x\notin \Q$. Then 
the set $\{\lceil mx\rceil-mx\}_{m\ge 1}$ is dense in
$[0,1]$ (cf.~\cite{Cassels57}, Chapter IV), hence 
$\sup_{m\in \Z}(\lceil mx\rceil-mx)=1$. We obtain
$1\le 1-x$, hence $x=0$. Contradiction.

Therefore $x=\frac{p}{q}$, for integers $1\le p\le q$ with 
$\gcd(p,q)=1$. The above inequality becomes
$$
1-\frac{1}{q}=\max_{m\in \Z}(\lceil mx\rceil-mx) \le 1-\frac{p}{q},
$$
hence $p=1$. Therefore $x=\frac{1}{q}$.
\end{proof}

We will need the following result of Lawrence.

\begin{thm}[\cite{Lawrence91}]\label{Law} 
Let $T=\R^d/\Z^d$ be a real torus.
\begin{itemize}
\item[(i)] Let $U\subset T$ be an open subset.
The the set of closed subgroups of $T$ which do not intersect $U$
has only finitely many maximal elements with respect to inclusion. 
\item[(ii)] The set of finite unions of closed subgroups of $T$
satisfies the descending chain condition.
\end{itemize}
\end{thm}

\begin{thm}\label{dioacc} 
Assume that $A$ satisfies the ascending chain condition. 
Then the set $\{\langle x,a\rangle; (x,a)\in \tilde{V}_d(A)\}$ 
satisfies the ascending chain condition.
\end{thm}

\begin{proof} Assume first that $d=1$. By Lemma~\ref{onedim},
$$
\{\langle x,a\rangle; (x,a)\in \tilde{V}_1(A)\}=
\{\frac{1}{n};n\ge 1\}\cdot A.
$$
Both sets $\{\frac{1}{n};n\ge 1\},A$ are nonnegative and 
satisfy the ascending chain condition, hence their product 
satisfies the ascending chain condition.

Let now $d\ge 2$ and assume by induction the result for smaller
values of $d$. Assume by contradiction that $\{(x^n,a^n)\}_{n\ge 1}$ 
is a sequence in $\tilde{V}_d(A)$ such that 
$$
\langle x^n,a^n\rangle<\langle x^{n+1},a^{n+1}\rangle
\mbox{ for } n\ge 1.
$$ 
Since $A$ satisfies the ascending chain condition, we may 
assume after passing to a subsequence that
$$
a_i^n\ge a_i^{n+1}, \forall n\ge 1, \forall 1\le i\le d.
$$
Assume first that $x^n\notin (0,1)^d$ for infinitely many $n$'s. 
After passing to a subsequence, we
may assume $x^n_1=1$ for every $n$. Write $x^n=(1,\bar{x}^n)$
and $a^n=(a_1^n,\bar{a}^n)$. Then 
$\langle \bar{x}^n,\bar{a}^n\rangle<\langle \bar{x}^{n+1},
\bar{a}^{n+1}\rangle$ for every $n\ge 1$, which contradicts
the acc property of the set
$\{\langle \bar{x},\bar{a}\rangle; (\bar{x},\bar{a})\in 
\tilde{V}_{d-1}(A)\}$.

Assume now that $x^n\in (0,1)^d$ for every $n$. We set 
$$
U^n=\{x\in (0,1)^d;\langle x-x^n,a^n\rangle<0\}
$$
and regard $U^n$ as an open subset of the torus 
$T^d=\R^d/\Z^d$. Let $X^n$ be the union of the subgroups of $T^d$ 
which do not intersect $U^n$. By Theorem~\ref{Law}.(i), $X^n$ is a 
finite union of closed subgroups of $T^d$. It is easy to see that 
$U^n\subseteq U^{n+1}$, hence $X^n\supseteq X^{n+1}$ for $n\ge 1$.

Since $(x^n,a^n)\in \tilde{V}_d(A)$, we have
$U^n\cap (\Z^d+\Z x^n)=\emptyset$. Therefore
$x^n\in X^n$ for every $n$. We have
$$
\langle x^n,a^{n+1}\rangle\le 
\langle x^n,a^n\rangle<\langle x^{n+1},a^{n+1}\rangle.
$$
Then $x^n\in U^{n+1}$, hence $x^n\notin X^{n+1}$. Therefore
$X^n\supsetneqq X^{n+1}$ for every $n\ge 1$, contradicting
Theorem~\ref{Law}.(ii).
\end{proof}

\begin{lem}\label{closure} 
The following properties hold:
\begin{itemize}
\item[(1)] If $A$ is a closed set, then $\tilde{V}_d(A)$ 
is a closed subset of $(0,1]^d\times A^d$.
\item[(2)] Identify $(0,1]^s$ with the face 
$x_{s+1}=\cdots=x_d=1$ of $(0,1]^d$. Then 
$$
\tilde{V}_d(A)\cap (0,1]^s=\tilde{V}_s(A).
$$
\item[(3)] Identify $[0,1]^s$ with the face 
$x_{s+1}=\cdots=x_d=0$ of $[0,1]^d$ and assume that $A$ 
is a closed set. Then 
$$
\overline{\tilde{V}_d(A)}\cap (0,1]^s=\tilde{V}_s(A).
$$
\end{itemize}
\end{lem}

\begin{proof} (1) Let $(x,a)\in (0,1]^d\times A^d$ such that
there exists a sequence $\{(x^n,a^n)\}_{n\ge 1}$ in 
$\tilde{V}_d(A)$ with $x=\lim_{n\to \infty} x^n$ and 
$a=\lim_{n\to \infty} a^n$. Fix $m\in \Z$. By assumption, 
we have
$$
\sum_{i=1}^d (1+(m-1)x_i^n-\lceil mx_i^n\rceil)a^n_i\ge 0, 
\forall n\ge 1.
$$
There exists a positive integer $n(m)$ such that
$
\lceil mx_i^n \rceil\ge \lceil mx_i\rceil
$ 
for every $1\le i\le d$ and every $n\ge n(m)$. Therefore
$$
\sum_{i=1}^d (1+(m-1)x_i^n-\lceil mx_i\rceil)a^n_i\ge 0, 
\forall n\ge n(m),
$$
Letting $n$ converge to infinity, we obtain 
$$
\sum_{i=1}^d (1+(m-1)x_i-\lceil mx_i\rceil)a_i\ge 0.
$$
Since $m$ was arbitrary, we conclude that $(x,a)\in \tilde{V}_d(A)$.

(2) This is clear.

(3) Assume that we have a sequence 
$\{(x^n,a^n)\}_{n\ge 1}\subset \tilde{V}_d(A)$ such that 
$\lim_{n\to \infty} x^n=(x,0,\ldots,0)\in (0,1]^s$ and 
$\lim_{n\to \infty}a^n=(a,a_{s+1},\ldots,a_d)$. 
Let $m$ be a positive integer. Note that for $s+1\le i\le d$ we 
have $mx^n_i\in (0,1]$ for $n\ge n(m)$, hence
$$
\lim_{n\to \infty} (1+(m-1)x_i^n-\lceil mx_i^n\rceil)=0
\mbox{ for } s+1\le i\le d.
$$
Therefore $\sum_{i=1}^s (1+(m-1)x_i-\lceil mx_i\rceil)a_i\ge 0$
for every $m\ge 1$. Since $\Z^s+\Z x$ is included
in the closure of $\Z^d+\Z_{\ge 0} x$, we obtain
$\sum_{i=1}^s (1+(m-1)x_i-\lceil mx_i\rceil)a_i\ge 0$ for
$m\le -1$ as well. Therefore $(x,a)\in \tilde{V}_s(A)$, 
proving the direct inclusion.

For the converse, just note that $(x,a)\in \tilde{V}_s(A)$ is 
the limit of the sequence 
$
((x,\frac{1}{n},\ldots,\frac{1}{n}),(a,1,\ldots,1))
\in \tilde{V}_d(A).
$
\end{proof}

\begin{defn} For $x\in \R$ and $m\in \Z$, define
$$
x^{(m)}=1+mx-\lceil mx\rceil.
$$
Note that this operation induces a selfmap of the 
half-open interval $(0,1]$. For $x\in \R^d$ and 
$m\in \Z$, define $x^{(m)}\in \R^d$ componentwise.
\end{defn}

Since $(0,1]^d\cap (\Z^d+\Z x)=\{x^{(m)}; m\in \Z\}$,
we have the equivalent description
$$
\tilde{V}_d(A)=\{(x,a)\in (0,1]^d\times A^d; 
\langle x^{(m)}-x,a\rangle\ge 0,\forall m\in \Z\}.
$$

\begin{lem}\label{limit} 
Let $x\in (0,1]^d$ and let $a\in A^d$ such that
$a_i>0$ for $1\le i\le d$. Then 
there exists a relatively open neighborhood $x\in U\subseteq 
(0,1]^d$ such that if $y\in U$ and
$\langle y^{(m)}-x,a\rangle\ge 0$ for every $m\in \Z$,
then $\langle y-x,a\rangle=0$.
\end{lem}

\begin{proof} (1) Assume first that $x\in (0,1)^d$.
By Theorem~\ref{Law}.(ii), the set of closed subgroups 
of $\R^d$ which contain $\Z^d$ and do not intersect 
the nonempty open set 
$
\{y\in (0,1)^d; \langle y-x,a\rangle <0\}
$
has finitely many maximal elements with respect to
inclusion, say $H_1,\ldots, H_l$.

If $x\in H_1$, then $H_1$ is a rational affine subspace of
$\R^d$ in an open neighborhood $x\in U_1\subset (0,1)^d$. 
Let $v\in H_1-x$. Since $x\in (0,1)^d$, there exists 
$\epsilon>0$ such that $x+tv\in H_1\cap (0,1)^d$ for 
$\vert t\vert<\epsilon$.
In particular, $\langle x+tv-x,a\rangle \ge 0$, that is
$t\langle v,a\rangle \ge 0$ for $\vert t\vert<\epsilon$.
We infer that $\langle v,a\rangle=0$. Therefore 
$H_1\cap U_1$ is contained in 
$\{y\in (0,1)^d; \langle y-x,a\rangle=0\}$.
If $x\notin H_1$, then $U_1=(0,1)^d\setminus H_1$ is an 
open neighborhood of $x$.

Repeating this procedure, we obtain a neighborhood $U_i$ 
of $x$, for each of the closed subgroup $H_i$. The 
intersection $U=U_1\cap \cdots \cap U_l$ is the desired 
neighborhood.

(2) We may assume after a reordering that 
$x_i=1$ for $1\le i\le s$ and $x_i\in (0,1)$ for 
$s< i\le n$. If $s=n$, we may take $U=(0,1]^d$. Assume 
now that $s<n$. By~\cite{Cassels57}, Chapter IV, there
exists a {\em negative} integer $m_0$ such that
$$
\langle x^{(m_0)}-x,a\rangle<\min_{i=1}^s a_i.
$$
Let $y\in (0,1]^s\times \prod_{i=s+1}^d 
(\frac{\lceil m_0x_i\rceil}{m_0}, \frac{\lceil m_0x_i\rceil-1}{m_0})$ 
such that $\langle y^{(m)}-x,a\rangle\ge 0$ for every $m\in \Z$.
We claim that $y_1=\cdots=y_s=1$. Indeed, assume by contradiction
that $y_j<1$ for some $1\le j\le s$.
A straightforward computation gives
$$
\langle y^{(m_0)}-x,a\rangle-m_0\langle y-x,a\rangle=
\langle x^{(m_0)}-x,a\rangle+\sum_{i=1}^d 
(\lceil m_0x_i\rceil-\lceil m_0y_i\rceil)a_i.
$$
By the choice of $y$, we obtain
$$
\sum_{i=1}^d (\lceil m_0x_i\rceil-\lceil m_0y_i\rceil)a_i
=\sum_{i=1}^s (m_0-\lceil m_0y_i\rceil)a_i\le -a_j,
$$
hence $0\le \langle x^{(m_0)}-x,a\rangle-a_j$. This 
contradicts our choice of $m_0$.

Let $\bar{x}=(x_{s+1},\ldots,x_d), \bar{y}=(y_{s+1},\ldots,
y_d), \bar{a}=(a_{s+1},\ldots,a_d)$. We have $(\bar{x},\bar{a})
\in \tilde{V}_{d-s}(A)$ and 
$\langle \bar{y}^{(m)}-\bar{x},\bar{a}\rangle\ge 0$ 
for every $m\in \Z$. From Step 1, there exists an open 
neighborhood $\bar{x}\in \bar{U}\subset (0,1)^s$ such that
if $\bar{y}\in \bar{U}$ then 
$\langle \bar{y}-\bar{x},\bar{a}\rangle=0$. Then 
$$
U=(0,1]^s\times (\bar{U}\cap \prod_{i=s+1}^d 
(\frac{\lceil m_0 x_i\rceil}{m_0}, \frac{\lceil m_0 x_i\rceil-1}{m_0})).
$$
satisfies the required properties.
\end{proof}

\begin{lem}\label{rat} 
The following equality holds for $d\ge 1$ and $a\in A^d$:
$$
\{\langle x,a\rangle;(x,a)\in \tilde{V}_d(A),x\in \Q^d\}=
\{\langle x,a\rangle; (x,a)\in \tilde{V}_d(A)\}.
$$
\end{lem}

\begin{proof} Let $(x,a)\in \tilde{V}_d(A)$. We have 
$x_1,\ldots,x_s<1$ and $x_{s+1}=\ldots=x_d=1$, where 
$0\le s\le d$. If $s=0$, then $x\in \Q^d$ and we are done.
Assume $s\ge 1$ and set $\bar{x}=(x_1,\ldots,x_s)$ and 
$\bar{a}=(a_1,\ldots,a_s)$. Then $(\bar{x},\bar{a})\in 
\tilde{V}_s(A)$. Since $\bar{x}\in (0,1)^s$, and there 
exists a closed subgroup $\Z^s\subseteq \bar{H}
\subseteq \R^s$ such that $\bar{x}\in \bar{H}
\cap U_{\bar{x}}\subset \{\bar{z}; \langle \bar{z}-\bar{x},
\bar{a}\rangle=0\}$, by Step 1 of the proof of 
Lemma~\ref{limit}.
Since $\bar{H}$ is rational, there exists $\bar{z}
\in \Q^s\cap H\cap U_{\bar{x}}$. Set $x'=(\bar{z},1,
\ldots,1)$. Then $(x',a)\in \tilde{V}_d(A)$,
$\langle x,a\rangle=\langle x',a\rangle$ and $x'\in \Q^d$.
\end{proof}

\begin{prop}\label{accum} 
Assume that $A$ has no positive accumulation points.
Then the set of accumulation points of 
$\{\langle x,a\rangle; (x,a)\in \tilde{V}_d(A)\}$
is 
$$
\{0\}\cup\bigcup_{1\le d'\le d-1}
\{\langle x,a\rangle; (x,a)\in \tilde{V}_{d'}(A)\}.
$$
\end{prop}

\begin{proof} Let $r>0$ be an accumulation point, that 
is there exists a sequence 
$(x^n,a^n)\in \tilde{V}_d(A)$ with 
$r=\lim_{n\to \infty}\langle x^n,a^n\rangle$ and
$r\ne \langle x^n,a^n\rangle$ for every $n\ge 1$.
By compactness, we may assume after passing to a subsequence
that $\lim_{n\to \infty} x^n=x\in [0,1]^d$ and 
$\lim_{n\to \infty} a^n=a\in [0,1]^d$ exist.
We have $r=\langle x,a\rangle$.

We claim that $a_ix_i=0$ for some $i$. Indeed, assume
by contradiction that $a_ix_i>0$ for every $1\le i\le d$.
Since $A$ has no nonzero accumulation points, we obtain 
$a^n=a$ for $n\ge 1$. Let $U_x\subset (0,1]^d$ be the 
relative neighborhood of $x$ associated to $(x,a)$ in 
Lemma~\ref{limit}. Then $x^n\in U_x$ for $n\ge n_0$. 
If $\langle x^n-x,a\rangle\ge 0$, then 
$(x^n,a)\in\tilde{V}_d(A)$ implies that 
$\langle z-x,a\rangle\ge 0$ for every 
$z\in (\Z^d+\Z x^n)\cap (0,1]^d$. Therefore
$\langle x^n-x,a\rangle=0$. This means 
$\langle x^n,a\rangle=r$, a contradiction.

Therefore $\langle x^n,a\rangle<r$ for every $n$.
Since $A$ has no positive accumulation points, it
satisfies the ascending chain condition. Therefore
the sequence $(\langle x^n,a\rangle)_{n\ge 1}$ satisfies
the ascending chain condition as well, by 
Theorem~\ref{dioacc}. This is a contradiction.

We may assume $a_ix_i>0$ for $1\le i\le {d'}$ and
$a_i x_i=0$ for $d'+1\le i\le d$. We have $d'\ge 1$,
since $\langle x,a\rangle>0$. Denote 
$\bar{x}=(x_1,\ldots,x_{d'})$ and 
$\bar{a}=(a_1,\ldots,a_{d'})$. We have 
$r=\langle \bar{x},\bar{a}\rangle$ and 
$(\bar{x},\bar{a})\in \tilde{V}_{d'}(A)$ by 
Theorem~\ref{closure}.

For the converse, note that
$$
((\frac{1}{k},\ldots,\frac{1}{k}), (1,\ldots,1))\in 
\tilde{V}_d(A)
$$
and $\langle (\frac{1}{k},\ldots,\frac{1}{k}), (1,\ldots,1)
\rangle=\frac{d}{k}$ accumulates to $0$. Let now 
$(x,a)\in \tilde{V}_{d'}(A)$ for $1\le d'\le d-1$.
Define $x_k=(x',\frac{1}{k},\ldots,\frac{1}{k})$ and
$a=(a',1,\ldots,1)$. Then $(x_k,a)\in \tilde{V}_d(A)$
and $\langle x_k,a\rangle=\langle x',a'\rangle+
\frac{d-d'}{k}$ accumulates to $\langle x',a'\rangle$.
\end{proof}

\begin{rem} Proposition~\ref{accum} is false if 
$A$ has a positive accumulation point. For example,
let $a>0$ be an accumulation point of a sequence of
elements $a_k\in A$.
Then $((1,\ldots,1),(a_k,1,\ldots,1))\in V_d(A)$
and $\langle (1,\ldots,1),(a_k,1,\ldots,1)\rangle$
accumulates to $d-1+a>d-1$, which clearly does not 
correspond to any element of $\tilde{V}_{d'}(A)$, for
$d'\le d-1$.
\end{rem}

There are many rational points in the set
$\tilde{V}_d\setminus V_d$. For example, $(\frac{1}{2},1)$ 
or $(\frac{l-1}{2l},\frac{1}{l})$ for $l\ge 2$. However,
the following property holds.

\begin{lem}\label{bartov} 
The following inclusion holds:
$$
\{\langle x,a\rangle;(x,a)\in \tilde{V}_d(A)\}
\subseteq 
\{\langle x,a\rangle; (x,a)\in V_d(\{\frac{1}{n};n\ge 1\}
\cdot A)\}.
$$
\end{lem}

\begin{proof} Let $r=\langle x,a\rangle$ for some 
$(x,a)\in \tilde{V}_d(A)$. By Lemma~\ref{rat}, we may assume
that $x\in \Q^d$. We may assume $a_i>0$ for every $i$.
Let $e_1,\ldots,e_d$ be the standard basis of $\R^d$,
spanning the standard cone $\sigma$, let 
$e=\sum_{i=1}^dx_ie_i$ and let $N=\sum_{i=1}^d \Z e_i+\Z e$.
If we set $\psi=\sum_{i=1}^d a_ie_i^*$, then we have
$$
\min(\psi\vert_{N\cap \relint(\sigma)})=\psi(e)=r.
$$ 
There exists positive integers $n_i\ge 1$ such that 
$e'_i=\frac{1}{n_i}e_i$ are primitive elements of the 
lattice $N$. In the new coordinates, we have
$\psi=\sum_{i=1}^d \frac{a_i}{n_i}{e'_i}^*$ and
$e=\sum_{i=1}n_ix_ie'_i$. Since $\psi$ attains its 
minimum at $e$ and all $a_i$'s are positive, we infer
that $n_ix_i<1$ for every $i$. Set 
$a'_i=\frac{a_i}{n_i}$ and $x'_i=n_ix_i$. Then
$(x',a')\in V_d(\{\frac{1}{n};n\ge 1\}\cdot A)$ and
$\langle x',a'\rangle=r$.
\end{proof}

\begin{cor} Assume that $A=\{\frac{1}{n};n\ge 1\}
\cdot A$. Then 
$$
\{\langle x,a\rangle;(x,a)\in V_d(A)\}
=
\{\langle x,a\rangle; (x,a)\in \tilde{V}_d(A)\}.
$$
\end{cor}
%%%%%%%%%%%%%%%%%%%%%%%%%%%%%%%%%%
%%%%%%%%%%%%%%%%%%%%%%%%%%%%%%%%%%

\section{Accumulation points of $\Mld^{tor}_d(A)$}

%%%%%%%%%%%%%%%%%%%%%%%%%%%%%%%%%%%
%%%%%%%%%%%%%%%%%%%%%%%%%%%%%%%%%%%

\begin{thm}\label{mt} The following properties hold:
\begin{itemize}
\item[(1)] If $A$ satisfies the ascending chain condition, 
then so does $\Mld^{tor}_d(A)$.
\item[(2)] Assume that $A$ has no positive accumulation 
points. Then the set of accumulation points of $\Mld^{tor}_d(A)$ 
is included in
$$
\{0\}\cup
\bigcup_{1\le d'\le d-1}\Mld^{tor}_{d'}(\{\frac{1}{n};n\ge 1\}\cdot A).
$$ 
The inclusion is an equality
if $\{\frac{1}{n};n\ge 1\}\cdot A\subset A$.
\item[(3)] Assume that $A$ has no positive accumulation
points and $\{\frac{1}{n};n\ge 1\}\cdot A\subset A$. 
Then $\Mld^{tor}_d(A)$ is a closed set if and only if 
$0\in A$.
\end{itemize}
\end{thm}

\begin{proof} (1) By Proposition~\ref{gutaiteki}, the 
inclusion $V_d(A)\subseteq \tilde{V}_d(A)$, and 
Theorem~\ref{dioacc}, the set $\Mld^{tor}_d(A)$ is a subset of 
a finite union of sets satisfying the ascending chain condition. 
Therefore $\Mld^{tor}_d(A)$ satisfies the ascending chain condition.

(2) Assume that $A$ has no nonzero accumulation points. 
By Proposition~\ref{accum} and Lemma~\ref{bartov}, the 
accumulation points of $\Mld^{tor}_d(A)$ belong to the set 
$$
\{0\}\cup
\bigcup_{1\le d'\le d-1}\Mld^{tor}_{d'}(\{\frac{1}{n};n\ge 1\}\cdot A).
$$
Assuming moreover that $\{\frac{1}{n};n\ge 1\}A\subset A$,
we will show that all points of the above set are accumulation
points of $\Mld^{tor}_d(A)$. If $(x,a)\in V_{d'}(A)$, 
then 
$$
((x,1,\ldots,1),(a,\frac{1}{n},\ldots,\frac{1}{n}))\in V_d(A),
$$
and $\langle (x,1,\ldots,1),
(a,\frac{1}{n},\ldots,\frac{1}{n})\rangle=\langle x,a\rangle+
\frac{d-d'}{n}$ accumulates to $\langle x,a\rangle$. Similarly,
$$
((1,1,\ldots,1),(\frac{1}{n},\frac{1}{n},\ldots,\frac{1}{n}))\in V_d(A)
$$ 
and
$\langle (1,1,\ldots,1),
(\frac{1}{n},\frac{1}{n},\ldots,\frac{1}{n})\rangle=\frac{d}{n}$ 
accumulates to $0$.
This proves the claim.

(3) Assume that $\Mld^{tor}_d(A)$ is a closed set. Since
$$
((\frac{1}{k},\ldots,\frac{1}{k}),(a,\ldots,a))\in V_d(A),
$$
we infer that $0=\lim_{k\to \infty}\frac{da}{k}\in \Mld^{tor}_d(A)$, 
which implies $0\in A$.

Conversely, assume $0\in A$. If $(x,a)\in V_{d'}(A)$
then 
$$
((x,1,\ldots,1),(a,0,\ldots,0))\in V_d(A)
$$
and $\langle (x,1,\ldots,1),(a,0,\ldots,0)\rangle =
\langle x,a\rangle$. We infer from (3) that $\Mld^{tor}_d(A)$ 
is a closed set.
\end{proof}

\begin{lem} Assume that $A$ has no positive accumulation
points. Then the following properties hold:
\begin{itemize}
\item[(1)] The set of accumulation points of $\Mld^{tor}_2(A)$ is 
$
\{0\}\cup \bigcup_{k\ge 1}\frac{1}{k}A.
$
\item[(2)] The set $\Mld^{tor}_2(A)$ is closed if and only
if $0\in A$.
\end{itemize}
\end{lem}

\begin{proof} (1) From Theorem~\ref{mt}, all accumulation
points are of this form. Conversely, fix $a\in A$ and 
$k\in \Z_{\ge 1}$. Then
$
((\frac{1}{kn+1},\frac{n}{nk+1}),(a,a))\in V_2(A)
$
is a sequence converging to $((0,\frac{1}{k}),(a,a))$,
hence $\frac{a}{k}$ is an accumulation point of
$\Mld^{tor}_2(A)$. Since $k$ is arbitrary, we infer
that $0$ is an accumulation point as well.

(2) Assume that $\Mld^{tor}_2(A)$ is a closed set.
Then $0\in \Mld^{tor}_2(A)$, which implies $0\in A$.

Assume now that $0\in A$. Then for $a\in A$ and 
$k\in \Z_{\ge 1}$ we have
$$
\frac{a}{k}=\langle (\frac{1}{k},\frac{1}{k}),(0,a)\rangle
\in \Mld^{tor}_2(A).
$$
From (1), these are all possible accumulation points, hence
$\Mld^{tor}_2(A)$ is a closed set.
\end{proof}

%%%%%%%%%%%%%%%%%%%%%%%%%%%%%%%%%%%%%%%
%%%%%%%%%%%%%%%%%%%%%%%%%%%%%%%%%%%%%%%

\end{document}